\documentclass[12pt]{article}
\begin{document}

\begin{center}
{\bf \LARGE A New Algorithm for the Higher-Order $G$-Transformation}
\end{center}
\vspace{2cm}

\begin{center}
{Avram Sidi\\
Computer Science Department\\
Technion - Israel Institute of Technology\\ Haifa 32000, Israel\\
\bigskip
e-mail:  asidi@cs.technion.ac.il\\
\vspace{1cm}
December 2000}
\end{center}
\newpage
\begin{abstract}
Let the scalars
$A^{(j)}_n$ be defined via the linear equations
$$A_l=A^{(j)}_n+\sum^n_{k=1}\bar{\alpha}_ku_{k+l-1},\ \
l=j,j+1,\ldots,j+n\ .$$ Here the $A_i$ and $u_i$ are known and the
$\bar{\alpha}_k$ are additional unknowns, and the quantities of
interest are the $A^{(j)}_n$. This problem arises, for example,
when one computes infinite-range integrals by the higher-order
$G$-transformation of Gray, Atchison, and McWilliams. One
efficient procedure for computing the $A^{(j)}_n$ is the
rs-algorithm of Pye and Atchison. In the present work,  we develop
yet another procedure that combines the FS-algorithm of Ford and
Sidi and the qd-algorithm of Rutishauser, and we denote it the
FS/qd-algorithm.  We show that the FS/qd-algorithm has a smaller
operation count than the rs-algorithm. We also show that the FS/qd
algorithm can also be used to implement the transformation of
Shanks, and compares very favorably with the
$\varepsilon$-algorithm of Wynn that is normally used for this
purpose.
\end{abstract}
\vspace{5cm}
\newpage
\pagenumbering{arabic}

\section{The Higher Order $G$-Transformation} \label{se1}
The $G$-transformation was designed by Gray and Atchison
\cite{Gray:1967:NTR} as an extrapolation method for
evaluating infinite
integrals of the form $\int^{\infty}_a f(t)\,dt\equiv I[f]$.
It was later
generalized in different ways in Atchison and Gray
\cite{Atchison:1968:NTR} and Gray and Atchison
\cite{Gray:1968:GT}, the ultimate generalization being given in
Gray, Atchison, and McWilliams \cite{Gray:1971:HOT}. This generalization
was denoted the higher-order $G$-transformation.
The way it is defined in \cite{Gray:1971:HOT},
this transformation produces
approximations to $I[f]$ that are of the form
\begin{equation} \label{eq:1}
G_n(x;h)= \frac{\left|
\begin{array}{cccc}
F(x)&F(x+h)&\cdots&F(x+nh)\\
f(x)&f(x+h)&\cdots&f(x+nh)\\
\vdots&\vdots&&\vdots \\
f(x+(n-1)h)&f(x+nh)&\cdots&f(x+(2n-1)h)
\end{array}
\right|}{\left|\begin{array}{cccc}
1&1&\cdots&1\\
f(x)&f(x+h)&\cdots&f(x+nh)\\
\vdots&\vdots&&\vdots \\
f(x+(n-1)h)&f(x+nh)&\cdots&f(x+(2n-1)h)
\end{array}
\right|},
\end{equation}
where
\begin{equation} \label{eq:2}
F(x)=\int^x_af(t)\, dt.
\end{equation}
The approximations produced by the $G$-transformation of
\cite{Gray:1967:NTR} are simply the $G_1(x;h)$.

It follows from (\ref{eq:1}) that $G_n(x;h)$ is also the
solution of the linear system
\begin{equation} \label{eq:3}
F(x+ih)=G_n(x;h)+\sum^n_{k=1}\bar{\alpha}_k f(x+(i+k-1)h),
\ \ i=0,1,\ldots, n,
\end{equation}
where $\bar{\alpha}_k$ are additional unknowns.

It has been shown in \cite{Gray:1971:HOT} that the kernel of the
higher-order $G$-transformation is the set of functions $f(x)$ that
are integrable at infinity in the sense of Abel and that
satisfy linear homogeneous ordinary differential equations of
order $n$ with constant coefficients. Thus $f(x)$ is in this kernel
if it is of the form $f(x)=\sum^r_{k=1}p_k(x)e^{c_k x}$, where
the  $c_k\neq 0$ are distinct and $\Re c_k\leq 0$,
$p_k(x)$ are polynomials.
If $\mu_k$ is the degree of $p_k(x)$ for each $k$, and if
$\sum^r_{k=1}(\mu_k+1)=n$, then $G_n(x;h)=I[f]$ for all $x$ and $h$.
On the basis of this result it was concluded in
Levin and Sidi \cite{Levin:1981:TNC}
that the higher-order $G$-transformation will be effective on functions
of the form $f(x)=\sum^s_{k=1}e^{c_kx}v_k(x)$, where $v_k(x)\sim
\sum^{\infty}_{i=0}\alpha_{ki}x^{\gamma_k-i}$ as $x\to\infty$,
with arbitrary $\gamma_k$, provided $h$ is
of a suitable size.

In the present work,  we are concerned with the actual computation
of the $G_n(x;h)$. Of course, it is not desirable to compute
$G_n(x;h)$ via the determinantal representation in (\ref{eq:1}).
Direct solution of the linear systems in (\ref{eq:3}) is expensive
too.  A very efficient and elegant procedure for computing the
$G_n(x;h)$ was given by Pye and  Atchison in \cite{Pye:1973:ACH},
and it has been denoted the rs-algorithm in Brezinski and Redivo
Zaglia \cite{Brezinski:1991:EMT}. The derivation of the
rs-algorithm makes use of the representation in (\ref{eq:1}). In
the present work,  we develop yet another procedure for computing
the $G_n(x;h)$, and we call this new procedure  the
FS/qd-algorithm. We show that the  FS/qd-algorithm is more
efficient than the rs-algorithm. With proper substitutions to be
discussed in Section \ref{se3}, the FS/qd-algorithm  can also be
used to implement the transformation of Shanks
\cite{Shanks:1955:NTD}, which is normally implemented via the
well-known $\varepsilon$-algorithm  of Wynn \cite{Wynn:1956:DCT}.
We show that, when used for implementing the transformation of
Shanks, the FS/qd-algorithm compares very favorably with the
 $\varepsilon$-algorithm.

\section{Algorithms for the Higher Order $G$-Transformation}
\label{se2}
We start with a review of the algorithm of
Pye and  Atchison
\cite{Pye:1973:ACH}.
Actually these authors consider the more general problem in which
one would like to compute the quantities $A^{(j)}_n$ defined via
the linear equations
\begin{equation} \label{eq:4}
 A_l=A^{(j)}_n+\sum^n_{k=1}\bar{\alpha}_ku_{k+l-1},\ \ l=j,j+1,\ldots,
j+n,
\end{equation}
where the $A_i$ and $u_i$ are known scalars, while the $\bar{\alpha}_k$
are not necessarily known.
Comparing (\ref{eq:3}) with (\ref{eq:4}) we can draw the analogy
$A_i\leftrightarrow F(x+ih)$, $u_i\leftrightarrow f(x+ih)$,
and $A^{(j)}_n\leftrightarrow G_n(x+jh;h)$.

The algorithm of \cite{Pye:1973:ACH} computes the $A^{(j)}_n$
with the help of
two sets of auxiliary quantities, $r^{(j)}_n$ and $s^{(j)}_n$.
These quantities are defined by
\begin{equation} \label{eq:5}
r^{(j)}_n=\frac{H^{(j)}_n}{K^{(j)}_n},\ \
s^{(j)}_n=\frac{K^{(j)}_{n+1}}{H^{(j)}_n},
\end{equation}
where $H^{(j)}_n$, the Hankel determinant
associated with $\{u_s\}$,
and $K^{(j)}_n$ are given as in
\begin{equation} \label{eq:6}
 H^{(j)}_n=\left|\begin{array}{cccc}
u_{j}&u_{j+1}&\cdots&u_{j+n-1}\\
u_{j+1}&u_{j+2}&\cdots&u_{j+n}\\
\vdots&\vdots &&\vdots\\
u_{j+n-1}&u_{j+n}&\cdots&u_{j+2n-2} \end{array}\right|,\ \
H^{(j)}_0=1,
\end{equation}
and
\begin{equation} \label{eq:7}
K^{(j)}_n=\left|\begin{array}{cccc}
1&1&\cdots&1\\
u_{j}&u_{j+1}&\cdots&u_{j+n-1}\\
u_{j+1}&u_{j+2}&\cdots&u_{j+n}\\
\vdots&\vdots &&\vdots\\
u_{j+n-2}&u_{j+n-1}&\cdots&u_{j+2n-3} \end{array}\right|,
\ \ K^{(j)}_0=1.
\end{equation}
The rs-algorithm computes the
$r^{(j)}_n$, and $s^{(j)}_n$ simultaneously by efficient recursions.
Once these have been computed the $A^{(j)}_n$ can be computed via
a separate recursion. \\

\noindent {\bf The rs-Algorithm} \hfill
\begin{enumerate}
\item
Set
$$s^{(j)}_0=1,\ \ r^{(j)}_1=u_j,\ \ A^{(j)}_0=A_j,\ \ j=0,1,\ldots\ .$$
\item
For $j=0,1,\ldots,$ and $n=1,2,\ldots,$ compute  recursively
$$s^{(j)}_{n}=s^{(j+1)}_{n-1}\left(\frac{r^{(j+1)}_{n}}{r^{(j)}_{n}}-1
\right), \ \ r^{(j)}_{n+1}=r^{(j+1)}_n\left(\frac{s^{(j+1)}_{n}}
{s^{(j)}_{n}}-1\right).$$
\item
For $j=0,1,\ldots,$ and $n=1,2,\ldots,$ set
$$A^{(j)}_n=\frac{r^{(j)}_{n}A^{(j+1)}_{n-1}-r^{(j+1)}_nA^{(j)}_{n-1}}
{r^{(j)}_{n}-r^{(j+1)}_n}. $$
\end{enumerate}

Before proceeding further, we note that
the equations in (\ref{eq:4}) are the same as
\begin{equation} \label{eq:8}
 A_l=A^{(j)}_n+\sum^n_{k=1}\bar{\alpha}_kg_k(l),\ \ l=j,j+1,\ldots,
j+n,
\end{equation}
with $g_k(l)=u_{k+l-1}$, $l=0,1,\ldots\ .$
Linear systems as in (\ref{eq:8})
arise also in the definition of a generalized
Richardson extrapolation method.
This suggests that the E-algorithm of Schneider
\cite{Schneider:1975:VRR} and the FS-algorithm of Ford and Sidi
\cite{Ford:1987:AGR} can be used for computing the $A^{(j)}_n$.
Different derivations of the E-algorithm were given by H{\aa}vie
\cite{Haavie:1979:GNT} and Brezinski \cite{Brezinski:1980:GEA}.
Of course, direct application of these algorithms
without taking into account the special nature of the $g_k(l)$ is
very uneconomical. By taking the nature of the $g_k(l)$ into
consideration, it becomes possible to derive  fast algorithms for
the $A^{(j)}_n$.

Now the E-algorithm produces the $A^{(j)}_n$ by a recursion relation
of the form
$$A^{(j)}_n=\frac{R^{(j)}_{n}A^{(j+1)}_{n-1}-R^{(j+1)}_nA^{(j)}_{n-1}}
{R^{(j)}_{n}-R^{(j+1)}_n}, $$
the most expensive part of the algorithm being the determination of
the $R^{(j)}_{n}$. Comparing the known expression for $R^{(j)}_{n}$
when $g_k(l)=u_{k+l-1}$ with that of $r^{(j)}_{n}$,
we realize that
$R^{(j)}_n=(-1)^{n-1}r^{(j)}_{n}$.
We thus conclude that the rs-algorithm is simply
the E-algorithm in which
the $R^{(j)}_{n}$, whose  determination forms the most
expensive part of the E-algorithm,  are computed by  a fast recursion.
For this point and others, see
\cite[Section 2.4]{Brezinski:1991:EMT}.

In view of the close connection between the rs- and E-algorithms, it is
natural to investigate the possibility of designing another algorithm
that is related to the FS-algorithm. This is worth the effort as
the FS-algorithm is more economical than the E-algorithm to begin with.
To this end we start with a brief description of the
FS-algorithm and refer
the reader to \cite{Ford:1987:AGR} for details.
For a comprehensive summary, see also Sidi \cite{Sidi:1990:GRE}.

Let us first define the short-hand notation
\begin{equation} \label{eq:9}
|u_1(j)\, u_2(j)\,\cdots\, u_n(j)|=\left |\begin{array}{cccc}
u_1(j)& u_2(j)&\cdots&u_n(j) \\ u_1(j+1)& u_2(j+1)&\cdots&u_n(j+1)\\
\vdots&\vdots&&\vdots\\
u_1(j+n-1)&u_2(j+n-1)&\cdots&u_n(j+n-1) \end{array}\right|,
\end{equation}
and set
\begin{equation} \label{eq:10}
G^{(j)}_n=|g_1(j)\, g_2(j)\,\cdots\, g_n(j)|, \ \ G^{(j)}_0=1.
\end{equation}
Next,  let us agree to denote the sequence $\{b(l)\}^{\infty}_{l=0}$
by $b$ for short, and define
\begin{equation} \label{eq:11}
f^{(j)}_n(b)=|g_1(j)\, g_2(j)\,\cdots\, g_n(j)\, b(j)|, \ \
f^{(j)}_0(b)=b(j).
\end{equation}
Finally, let us define
\begin{equation} \label{eq:12}
\psi^{(j)}_n(b)=\frac{f^{(j)}_n(b)}{G^{(j)}_{n+1}}.
\end{equation}
Then we have
\begin{equation} \label{eq:13}
A^{(j)}_n=\frac{\psi^{(j)}_n(a)}{\psi^{(j)}_n(I)},
\end{equation}
where $a$ and $I$ denote the sequences $\{a(l)=A_l\}^{\infty}_{l=0}$
and $\{I(l)=1\}^{\infty}_{l=0}$ respectively.
The FS-algorithm computes the quantities $\psi^{(j)}_n(b)$ by a
recursion of the form
\begin{equation} \label{eq:14}
\psi^{(j)}_n(b)=\frac{\psi^{(j+1)}_{n-1}(b)-\psi^{(j)}_{n-1}(b)}
{D^{(j)}_n}, \ \ D^{(j)}_n=\frac{G^{(j)}_{n+1}G^{(j+1)}_{n-1}}
{G^{(j)}_{n}G^{(j+1)}_{n}}.
\end{equation}
By this recursion,  first the $\psi^{(j)}_n(a)$ and $\psi^{(j)}_n(I)$ are
computed and then $A^{(j)}_n$ is determined via (\ref{eq:13}).

We recall that the most expensive part of the FS-algorithm
is the (recursive) determination of the quantities
$D^{(j)}_n$, and we would like to
reduce the  cost of this part.
Fortunately, this can be achieved once we
realize that, with $G^{(j)}_n$ and $H^{(j)}_n$ defined as in
(\ref{eq:10}) and (\ref{eq:6}) respectively,  and with
$g_k(l)=u_{k+l-1}$, we have
$G^{(j)}_n=H^{(j)}_n$ in the present case.
From this and (\ref{eq:14}),
we obtain the surprising result that

\begin{equation} \label{eq:15}
D^{(j)}_n=\frac{H^{(j)}_{n+1}H^{(j+1)}_{n-1}}
{H^{(j)}_{n}H^{(j+1)}_{n}}=e^{(j)}_n.
\end{equation}
Here $e^{(j)}_n$ is a quantity computed by the famous qd-algorithm of
Rutishauser \cite{Rutishauser:1954:QD1}, a clear exposition of which
can also be found in Henrici \cite{Henrici:1974:ACC}.
Actually, the qd-algorithm computes along with the $e^{(j)}_n$
also the quantities $q^{(j)}_n$ that  are given as in
\begin{equation} \label{eq:16}
q^{(j)}_n=\frac{H^{(j)}_{n-1}H^{(j+1)}_{n}}
{H^{(j)}_{n}H^{(j+1)}_{n-1}}.
\end{equation}
The $q^{(j)}_n$ and $e^{(j)}_n$ serve the construction of
the regular $C$-fractions, hence the Pad\'{e} approximants,  associated
with the formal power series $\sum^{\infty}_{k=0}u_kz^k$.
Regular $C$-fractions are continued fractions of a special type.
(See Jones and Thron \cite{Jones:1980:CFA}.)
The qd-algorithm computes the  $q^{(j)}_n$ and $e^{(j)}_n$
via the recursions
\begin{equation} \label{eq:17}
e^{(j)}_n=q^{(j+1)}_n-q^{(j)}_n+e^{(j+1)}_{n-1},\ \
q^{(j)}_{n+1}=\frac{e^{(j+1)}_{n}}{e^{(j)}_n}q^{(j+1)}_{n}, \ \
j\geq 0,\ n\geq 1,
\end{equation}
with the initial conditions $e^{(j)}_0=0$ and $q^{(j)}_1=u_{j+1}/u_j$
for all $j\geq 0$.
The quantities $q^{(k)}_n$ and $e^{(k)}_n$ can be arranged in a
two-dimensional array as in Figure \ref{fig}.
\begin{figure}
\[ \begin{array}{cccccc}
& q^{(0)}_1 & & & &\\
e^{(1)}_0  & & e^{(0)}_1\\
& q^{(1)}_1 & & q^{(0)}_2\\
e^{(2)}_0  & & e^{(1)}_1 & & e^{(0)}_2\\
& q^{(2)}_1 & & q^{(1)}_2 & &\ddots\\
e^{(3)}_0  & & e^{(2)}_1 & & e^{(1)}_2\\
& q^{(3)}_1 & & q^{(2)}_2 & & \ddots\\
e^{(4)}_0  & \vdots & e^{(3)}_1 & &  e^{(2)}_2\\
\vdots & & \vdots & q^{(3)}_2 & \vdots & \ddots\\
 & & & \vdots &
 \end{array} \]
\caption
{\label{fig} The qd-table.}
\end{figure}

This observation enables us to combine the FS- and qd-algorithms
to obtain the following economical implementation, the
FS/qd-algorithm, for the higher-order $G$-transformation.
For simplicity of notation, we will let $\psi^{(j)}_n(a)=M^{(j)}_n$
and $\psi^{(j)}_n(I)=N^{(j)}_n$ in the FS-algorithm. \\

\noindent {\bf The FS/qd-Algorithm}
\begin{enumerate}
\item
Set
$$e^{(j)}_0=0,\ \ q^{(j)}_1=\frac{u_{j+1}}{u_j},\ \
M^{(j)}_0=\frac{A_j}{u_{j}},\ \
N^{(j)}_0=\frac{1}{u_{j}}. $$
\item
For $j=0,1,\ldots,$ and $n=1,2,\ldots,$ compute recursively
$$e^{(j)}_n=q^{(j+1)}_n-q^{(j)}_n+e^{(j+1)}_{n-1},\ \
q^{(j)}_{n+1}=\frac{e^{(j+1)}_{n}}{e^{(j)}_n}q^{(j+1)}_{n},$$
$$M^{(j)}_n=\frac{M^{(j+1)}_{n-1}-M^{(j)}_{n-1}}
{e^{(j)}_n}, \ \
N^{(j)}_n=\frac{N^{(j+1)}_{n-1}-N^{(j)}_{n-1}}
{e^{(j)}_n}.$$
\item
For $j=0,1,\ldots,$ and $n=1,2,\ldots,$ set
$$A^{(j)}_n=\frac{M^{(j)}_n}{N^{(j)}_n}.$$
\end{enumerate}

It seems that, when a certain number of the $u_i$ are given,
it is best to compute the associated qd-table columnwise.
(Note that  the quantities in each of the recursions
for $e^{(j)}_n$ and $q^{(j)}_{n+1}$ in step 2 of the FS/qd-algorithm
above form the four corners of a lozenge
in Figure \ref{fig}.)
Following that we can compute the $M^{(j)}_n$, $N^{(j)}_n$,
and $A^{(j)}_n$ columnwise as well.

\section{Comparison of the rs- and FS/qd-Algorithms}\label{se3}
Let us now compare the operation counts of the two algorithms.
First, we note that the $r^{(j)}_n$ and $s^{(j)}_n$ in the
rs-algorithm can be arranged in
a table similar to the qd-table of the $e^{(j)}_n$ and $q^{(j)}_n$.
Thus, given $A_0,A_1,\ldots,A_L$, and $u_0,u_1,\ldots,u_{2L}$,
we can compute $A^{(j)}_n$ for $0\leq j+n\leq L$.
Now the number of the $e^{(j)}_n$ in the relevant qd-table is
$L^2+O(L)$ and so is that of the $q^{(j)}_n$.
A similar statement can be made concerning the $r^{(j)}_n$ and
$s^{(j)}_n$. The number of the $A^{(j)}_n$ is $L^2/2+O(L)$, and so
are the numbers of the $M^{(j)}_n$ and the $N^{(j)}_n$.
Consequently, we have the following operation counts.
\vspace{4ex}

\begin{tabular}{llll}
Algorithm & No. of Multiplications & No. of Additions & No. of Divisions
\\ \hline
FS/qd & $L^2+O(L)$ & $3L^2+O(L)$ & $5L^2/2+O(L)$\\
rs & $3L^2+O(L)$ & $3L^2+O(L)$ & $5L^2/2+O(L)$
\end{tabular}
\vspace{2ex}

\noindent In case only the $A^{(0)}_n$  are needed (as they have the best
convergence properties), the number of divisions in the FS/qd-algorithm
can be reduced from $5L^2/2+O(L)$ to $2L^2+O(L)$. In any case,  we see
that the operation count of the rs-algorithm  is about 30\% more
than that of the FS/qd-algorithm.

Finally, we observe that when $u_k=\Delta A_k=A_{k+1}-A_k$
the higher-order $G$-transformation reduces to the Shanks
\cite{Shanks:1955:NTD}
transformation, and therefore
the rs- and FS/qd-algorithms can be used for
computing the $A^{(j)}_n$ in this case too.
Of course, the most famous and efficient implementation of the
Shanks transformation is via the $\varepsilon$-algorithm of Wynn
\cite{Wynn:1956:DCT}, which reads
\begin{equation} \label{eq:18}
\varepsilon^{(j)}_{k+1}=\varepsilon^{(j+1)}_{k-1}+\frac{1}
{\varepsilon^{(j+1)}_{k}-\varepsilon^{(j)}_{k}}, \ \
j,k=0,1,\ldots,
\end{equation}
with the initial values $\varepsilon^{(j)}_{-1}=0$ and
$\varepsilon^{(j)}_{0}=A_j$ for all $j\geq 0$.
Then $\varepsilon^{(j)}_{2n}=A^{(j)}_n$ for all $j$ and $n$.
Thus, given $A_0,A_1,\ldots,A_{2L}$,
we can compute $\varepsilon^{(j)}_{2n}$ for $0\leq j+2n\leq 2L$.
(In particular, we can compute the diagonal approximations
$\varepsilon^{(j)}_{2n}$, $n=0,1,\ldots,L,$ that have the best
convergence properties.)
Since there are $2L^2+O(L)$ of the $\varepsilon^{(j)}_k$ to compute,
the operation count of this computation is
$4L^2+O(L)$ additions and $2L^2+O(L)$ divisions and no multiplications.
It is seen from the table above that the FS/qd-algorithm compares
very favorably with the $\varepsilon$-algorithm as well.


\end{document}